\documentclass[11pt]{article}
\usepackage{amsmath}
\usepackage{amssymb}

\font\tenmsb=msbm10

\def\eps{\varepsilon}
\font\tencmmib=cmmib10 \skewchar\tencmmib '60
\newfam\cmmibfam
\textfont\cmmibfam=\tencmmib

\font\tenmsb=msbm10

\def\Bbb#1{\hbox{\tenmsb#1}}

\def\bbox{\quad\hbox{\vrule \vbox{\hrule \vskip2pt \hbox{\hskip2pt
\vbox{\hsize=1pt}\hskip2pt} \vskip2pt\hrule}\vrule}}
\def\lessim{\ \lower4pt\hbox{$
\buildrel{\displaystyle <}\over\sim$}\ }
\def\gessim{\ \lower4pt\hbox{$\buildrel{\displaystyle >}
\over\sim$}\ }

\def\eps{{\varepsilon}}

\def\Bbb E{\mathbb{E}}
\def\Bbb R{\mathbb{R}} 
\newtheorem{lemma}{Lemma}
\newtheorem{theorem}{Theorem}
\newtheorem{corollary}{Corollary}
\makeatletter
\@addtoreset{equation}{section}

\makeatother

\font\tencmmib=cmmib10 \skewchar\tencmmib '60
\newfam\cmmibfam
\textfont\cmmibfam=\tencmmib

\font\tenmsb=msbm10

\def\Bbb#1{\hbox{\tenmsb#1}}

\def\bbox{\quad\hbox{\vrule \vbox{\hrule \vskip2pt \hbox{\hskip2pt
\vbox{\hsize=1pt}\hskip2pt} \vskip2pt\hrule}\vrule}}
\def\lessim{\ \lower4pt\hbox{$
\buildrel{\displaystyle <}\over\sim$}\ }
\def\gessim{\ \lower4pt\hbox{$\buildrel{\displaystyle >}
\over\sim$}\ }

\def\EE{{\Bbb E} }

\def\eps{\varepsilon}

\def\go0{\to 0}

\def\leftitem#1{\item{\hbox to\parindent{\enspace#1\hfill}}}

\def\qed{{$\hfill \bbox$}}

\def\sg{\sigma}

\def\sg2{\sigma^2}

\def\__{_{\infty}}

\begin{document}
{\baselineskip=16.5pt

\title{\bf Low Rank Estimation of Similarities on Graphs}

\author{
{\bf Vladimir Koltchinskii}
\thanks{Partially supported by NSF Grants DMS-1207808, DMS-0906880 and CCF-0808863}
\ \ {\bf and Pedro Rangel}
\thanks{Supported by NSF Grant CCF-0808863}
\\ School of Mathematics
\\ Georgia Institute of Technology
\\ Atlanta, GA 30332-0160
\\ vlad@math.gatech.edu, prangel@math.gatech.edu
}

\maketitle

\begin{abstract}
Let $(V,E)$ be a graph with vertex set $V$ and edge set $E.$ 
Let $(X,X',Y)\in V\times V\times \{-1,1\}$ be 
a random triple, where $X,X'$ are independent uniformly distributed vertices 
and $Y$ is a label indicating whether $X,X'$ are ``similar'' ($Y=+1$), or not ($Y=-1$).  
Our goal is to estimate the regression function 
$$S_{\ast}(u,v)={\mathbb E}(Y|X=u,X'=v), u,v\in V$$ 
based on training data consisting of $n$ i.i.d.
copies of $(X,X',Y).$ We are interested in this problem in the case when $S_{\ast}$ is a symmetric low rank kernel and, in addition to this, it is assumed that $S_{\ast}$ is ``smooth'' on the graph. 
We study estimators based on a modified least squares method with complexity penalization involving both the nuclear
norm and Sobolev type norms of symmetric kernels on the graph and prove upper  
bounds on $L_2$-type errors of such estimators with explicit dependence both on the rank of $S_{\ast}$ and on the degree of its smoothness. 
\end{abstract}



\medskip

\section{Introduction}

\medskip


\medskip

Let $G=(V,E)$ be a graph with vertex set $V$ and edge set $E,$ ${\rm card}(V)=m.$
Let $A:=(a(u,v))_{u,v\in V}$ be the adjacency matrix of $G,$ that is, $a(u,v)=1$ 
if $u$ and $v$ are connected with an edge and $a(u,v)=0$ otherwise.
Let $\Delta:=D-A$ be the Laplacian of $G,$ $D$ being the diagonal matrix with the degrees of vertices on the diagonal. Let $(X,X',Y)\in V\times V\times \{-1,1\}$ be a random triple with $X,X'$ being independent vertices sampled at random from the uniform distribution $\Pi$ on $V$ and $Y$ being an ``indicator'' of a symmetric binary relationship between $X,X'$ called in what follows a ``similarity''. More precisely, $Y=+1$ indicates that the vertices $X,X'$ are similar and $Y=-1$ indicates that they are not. The conditional distribution of $Y$ given $X,X'$ is completely characterized by the regression function 
$$
S_{\ast}(u,v):= {\mathbb E}(Y|X=u,X'=v), u,v\in V
$$
that is assumed to be a symmetric kernel on $V\times V$ and will be called the \emph{similarity kernel}. It is well known that ${\rm sign}(S_{\ast}(X,X'))$ is the Bayes classifier, that is, the best possible predictor of $Y$ based on an observation of $X,X'$ in the sense that it minimizes the generalization error ${\mathbb P}\{Y\neq g(X,X')\}$ over all possible predictors $g:V\times V\mapsto \{-1,1\}.$ Our goal is to estimate $S_{\ast}$ based on the training data $(X_1,X_1',Y_1),\dots, (X_n,X_n',Y_n)$ consisting of $n$ i.i.d. copies of $(X,X',Y).$ We are especially interested in the class of problems such that, on the one hand,  $S_{\ast}$ is a matrix (kernel) of relatively small rank and, on the other hand, $S_{\ast}$ possesses certain degree of smoothness on the graph.

Throughout the paper, $\mathcal{S}_V$ denotes the linear space of \emph{symmetric kernels} $S:V\times V\mapsto {\mathbb R},$ $S(u,v)=S(v,u), u,v\in V,$ that can be 
also viewed as real-valued symmetric $m\times m$ matrices. For $S \in \mathcal{S}_V$, let ${\rm rank}(S)$ denote the \emph{rank} of $S$ and ${\rm tr}(S)$
denote the \emph{trace} of $S.$
The \emph{spectral representation} of $S$ has the form $S=\sum_{j=1}^{r}{\sigma_j (\psi_j \otimes \psi_j)},$ where $r={\rm rank}(S)$, $\sigma_1\leq \dots \leq \sigma_r$ are non-zero eigenvalues of $S$ (repeated with their multiplicities) and $\psi_1,\dots,\psi_r$ are the corresponding orthonormal eigenfunctions (there is a multiple choice of $\psi_j$s in the case of repeated eigenvalues). We also use the notation ${\rm sign}(S):=\sum_{j=1}^{r}{{\rm sign}(\sigma_j)(\psi_j \otimes \psi_j})$ and we define the \emph{support} of $S$, denoted by ${\rm supp}(S)$, as the linear span of $\{ \psi_1,\dots,\psi_r\}$ in ${\mathbb R}^V.$

For $1\leq p< \infty,$ the Schatten $p$-norm of $S\in {\mathcal S}_V$ is defined as
$$
\|S\|_p:=({\rm tr}(|S|^p))^{1/p}=\left( \sum_{j=1}^{r} |\sigma_j|^p \right)^{1/p},
$$
where $|S|:=\sqrt{S^2}.$
For $p=1,$ $\|\cdot\|_1$ is called the \emph{nuclear norm}, while, for $p=2,$ $\|\cdot\|_2$ is the \emph{Hilbert--Schmidt} or \emph{Frobenius norm}, that is, the norm induced by the Hilbert--Schmidt inner product which will be denoted by $\langle \cdot,\cdot\rangle .$ The \emph{operator} or \emph{spectral norm} is defined as  $\|S\|:=\max_{j}|\sigma_j|.$  

Let us also denote by $\Pi^2:=\Pi\times \Pi$ the distribution of random 
couple $(X,X')$ in $V\times V$ and let $\|S\|_{L_2(\Pi^2)}$ be the 
$L_2(\Pi^2)$-norm of kernel $S:$
$$
\|S\|_{L_2(\Pi^2)}^2=\int_{V\times V}|S(u,v)|^2\Pi^2(du,dv)={\mathbb E}|S(X,X')|^2.
$$
The corresponding inner product is denoted by 
$\langle \cdot,\cdot \rangle_{L_2(\Pi^2)}.$
Clearly, under the assumption that the distribution $\Pi$ is uniform in $V,$ we have 
$\|S\|_{L_2(\Pi^2)}=m^{-2}\|S\|_2^2$ and 
$\langle S_1,S_2 \rangle_{L_2(\Pi^2)}=m^{-2}\langle S_1,S_2\rangle.$



The smoothness of a symmetric kernel $S:V\times V\mapsto {\mathbb R}$ can be characterized in terms of Sobolev type norms $\|\Delta^{p/2}S\|_2^2$ for some $p>0.$ Note that if $S$ is a kernel of rank $r$ with spectral representation $S=\sum_{k=1}^r \mu_k (\psi_k \otimes \psi_k)$, then
\footnote{Below $\|\cdot\|$ denotes the Euclidean norm in ${\mathbb R}^V;$
there is a little abuse of notation here since we also denote the operator norm
by $\|\cdot\|.$} 
$$
\|\Delta^{p/2}S\|_2^2=
{\rm tr}(\Delta^{p/2}S^2\Delta^{p/2})={\rm tr}(\Delta^pS^2)=
\sum_{k=1}^m \mu_k^2 \langle \Delta^p \psi_k,\psi_k\rangle=
\sum_{k=1}^m \mu_k^2 \|\Delta^{p/2} \psi_k\|^2,
$$
so, essentially, the smoothness of the kernel $S$ depends on the smoothness of its eigenfunctions $\psi_k$ on the graph. In particular, for $p=1,$ we have 
$$
\|\Delta^{1/2}S\|_2^2=
\sum_{k=1}^m \mu_k^2 \sum_{u\sim v} |\psi_k(u)-\psi_k(v)|^2,
$$
where the sum is over the couples of vertices connected with an edge. 

Given a kernel $S$, let $L_n(S)$ denote the following penalized empirical risk: 
\begin{equation}\label{linearized_LS}
\begin{split}
L_n(S)&:=
\biggl[
\|S\|_{L_2(\Pi^2)}^2-\frac{2}{n}\sum_{j=1}^n Y_j S(X_j,X_j')+ 
\eps \|S\|_1 + \bar{\eps} \|W^{1/2}S\|_{L_2(\Pi^2)}^2
\biggr] \\
&=
\biggl[
\|S\|_{L_2(\Pi^2)}^2-\frac{2}{n}\sum_{j=1}^n Y_j S(X_j,X_j')+ 
\eps \|S\|_1 + \eps_1 \|W^{1/2}S\|_2^2
\biggr]
\end{split}
\end{equation}
where $W=d\Delta^p$ for some constants $d>0$ and $p>0$, $\eps, \bar{\eps}>0$ are regularization parameters and $\eps_1=\frac{\bar \eps}{m^2}$. We will study the following estimation method:
\begin{equation} 
\label{nuce_min}
\hat S:={\rm argmin}_{S \in \mathbb{D}} L_n(S),
\end{equation}
where ${\mathbb D}$ is a closed convex subset of the linear space ${\cal S}_V$ of all symmetric kernels. Note that there are two complexity penalties involved in the definition of penalized empirical risk (\ref{linearized_LS}). The first penalty is based on the nuclear norm $\|S\|_1$ and it is used to ``promote'' low rank solutions. The second penalty is based on a ``Sobolev type norm'' $\|W^{1/2}S\|_2^2.$ It is used to ``promote'' the smoothness of the solution on the graph. In principle, $W$ in the definition of $L_n(S)$ could be an arbitrary symmetric nonnegatively definite matrix. Therefore, alternative interpretations of the problem under consideration are possible (such as, for instance, learning 
similarities on weighted graphs).

We will derive an upper bound on the error 
$\|\hat S-S_{\ast}\|_{L_2(\Pi^2)}^2=m^{-2}\|\hat S-S_{\ast}\|_2^2$ of estimator $\hat S$ in terms of spectral characteristics of the target similarity matrix $S_{\ast}$ and matrix $W$. Before stating the main results, let us recall recent advances on low rank matrix completion problems in which the approach based on nuclear norm penalization has been crucial. 

Suppose first that a symmetric kernel $S_{\ast}\in {\cal S}_V$ is 
observed at random points $(X_j,X_j'), j=1,\dots, n,$ where $X_j,X_j',j=1,\dots,n$
are independent and sampled from the uniform distribution $\Pi$ in $V.$
In this case, $V$ is an arbitrary finite set of cardinality $m$ and 
the set of edges $E$ is not specified. It is assumed that $Y_j=S_{\ast}(X_j,X_j'),$
so, there is no errors in the observations. In such a noiseless case, the following 
method is used to recover $S_{\ast}$ based on the observations $(X_1,X_1',Y_1),\dots, (X_n,X_n',Y_n):$
$$
\check S:={\rm argmin}\{\|S\|_1: S\in {\cal S}_V, S(X_j,X_j')=Y_j, j=1,\dots, n\}.
$$
Such methods of recovery of low rank target matrices $S_{\ast}$ have been 
extensively studied in the recent literature (see Candes and Recht (2009), Recht, Fazel and Parrilo (2010), Candes and Tao (2010), Gross (2011) and references 
therein). It is easy to see that there are low rank matrices $S_{\ast}$ 
that can not be recovered based on a random sample of $n$ entries unless 
$n$ is very large (comparable with the total number of entries of the matrix).
Indeed, consider $S_{\ast}$ such that, for given $u,v\in V,$ $S_{\ast}(u,v)=S_{\ast}(v,u)=1$ and $S_{\ast}(u',v')=0$ otherwise. 
For this rank $2$ matrix, the probability that the two ``informative''
entries are not present in the sample is $(1-\frac{2}{m^2})^n,$ which is 
close to $1$ if $n=o(m^2).$ Such sparse low rank matrices should be excluded
to make it possible to recover the target low rank matrix based 
on relatively small samples of entries. This is done by introducing so called 
\emph{low coherence} assumptions. Let $\{e_v:v\in V\}$ be the canonical orthonormal
basis of ${\mathbb R}^V$ equipped with the standard Euclidean inner product.
Given a linear subspace $L\subset {\mathbb R}^V,$ denote by $L^\perp$ the orthogonal complement of $L$ and by $P_L$ the projector onto the subspace $L.$ Let $L:={\rm supp}(S_{\ast}),$ $r={\rm rank}(S_{\ast})$ and suppose there exists a constant $\nu\geq 1$ (\emph{coherence coefficient}) such that 
\begin{equation}
\label{coherA}
\|P_L e_v\|^2 \leq \frac{\nu r }{m},\ v\in V
\ \ {\rm and}\ \ 
|\langle {\rm sign}(S_{\ast})e_u, e_v\rangle|^2 \leq \frac{\nu r}{m^2}, u,v\in V.
\end{equation}

The following result is due to Candes and Tao (2010) and Gross (2011)
(we state here a version of Gross that is an improvement of an earlier 
result of Candes and Tao with significant simplification of the proof). 

\begin{theorem}
\label{gross}
Suppose conditions (\ref{coherA}) hold for some 
$\nu\geq 1.$ Then, there exists a constant $C>0$ such that, for all 
$n\geq C\nu r m \log^2 m,$ $\check S=S_{\ast}$ with probability at least $1-m^{-2}.$
\end{theorem}

Thus, if, for the target matrix $S_{\ast},$ the coherence coefficient $\nu\geq 1$ is relatively small, the nuclear norm minimization algorithm (\ref{nuce_min})
does provide the exact recovery of $S_{\ast}$ as soon as the number of observed
entries $n$ is of the order $mr$ (up to a log factor).

In the case when $Y_j$ are noisy observations of $S_{\ast}(X_j,X_j')$
with 
$${\mathbb E}(Y_j|X_j=u,X_j'=v)=S_{\ast}(u,v),$$
one can use the following 
estimation method based on penalized empirical risk minimization with quadratic 
loss and with nuclear norm penalty:
\begin{equation}
\label{LS_nuce}
\check S:={\rm argmin}_{S\in {\cal S}_V}\biggl[n^{-1}\sum_{j=1}^n (Y_j-S(X_j,X_j'))^2 +\eps \|S\|_1\biggr].
\end{equation}
This method has been also extensively studied for the recent years,
in particular, by Candes and Plan (2011), Rohde and Tsybakov (2011),
Negahban and Wainwright (2010),
Koltchinskii, Lounici and Tsybakov (2011), Koltchinskii (2011b).
It was also pointed out by Koltchinskii, Lounici and Tsybakov (2011)
that in the case of known design distribution $\Pi$ (which is the case 
in our paper) one can use instead of (\ref{LS_nuce}) the following 
modified method: \footnote{Note that, if the norm $\|S\|_{L_2(\Pi^2)}$ in the definition 
below is replaced by the $L_2(\Pi_n)$-norm, where $\Pi_n$ is the empirical 
distribution based on $(X_1,X_1'),\dots, (X_n,X_n'),$ then the resulting 
estimator coincides with (\ref{LS_nuce}).}
\begin{equation}
\label{LS_nuce'}
\check S:={\rm argmin}_{S\in {\cal S}_V}\biggl[\|S\|^2_{L_2(\Pi^2)} 
-\frac{2}{n}\sum_{j=1}^n Y_j S(X_j,X_j')+\eps \|S\|_1\biggr].
\end{equation}
Clearly, (\ref{LS_nuce'}) is equivalent to method (\ref{nuce_min}) 
defined above for $\bar\eps=0.$

When the observations $|Y_j|\leq 1,j=1,\dots,n$ (for instance, 
when $Y_j\in \{-1,1\},$ which is the case studied in the paper), the next
result follows from Theorem 4 in Koltchinskii, Lounici and Tsybakov (2011).

\begin{theorem}
\label{K-L-T}
For $t>0,$ suppose that 
$$
\eps \geq 4 \biggl(\sqrt{\frac{t+\log (2m)}{nm}}\bigvee \frac{2(t+\log (2m))}{n}\biggr).
$$
Then with probability at least $1-e^{-t}$
$$
\|\check S-S_{\ast}\|_{L_2(\Pi)}^2 \leq \Bigl(\frac{1+\sqrt{2}}{2}\Bigr)^2 m^2 \eps^2 {\rm rank}(S_{\ast}).
$$
\end{theorem}

Our main goal is to show that this bound can be improved in the case when
the target kernel $S_{\ast},$ in addition to having relatively small rank, is also smooth on the 
graph and when the estimation method (\ref{nuce_min}) is used with a proper choice 
of regularization parameters $\eps, \bar \eps.$

\section{Main Results}

Suppose that $W$ has the following spectral representation:
$
W=\sum_{k=1}^m \lambda_k (\phi_k \otimes \phi_k),
$
where $0\leq \lambda_1 \leq \dots \leq \lambda_m$ are the eigenvalues of $W$ (repeated with their multiplicities) and $\phi_1,\dots, \phi_m$ are the corresponding orthonormal eigenfunctions (of course, there is a multiple choice of $\phi_k$ in the case of repeated eigenvalues). Let $k_0$ be the smallest $k$ such that $\lambda_k> 0$. We will assume that for some (arbitrarily large) $\zeta\geq 1$
$\lambda_m \leq m^{\zeta}$ and $\lambda_{k_0}\geq m^{-\zeta}.$ 
In addition, it is assumed that $s\mapsto \frac{s}{\lambda_s}$ is a nonincreasing sequence, that, for all $k=k_0,\dots, m-1,$ $\lambda_{k+1}\leq c\lambda_k,$
and, that, for all $s\geq k_0,$ 
\begin{equation}
\label{spec}
\sum_{k=s}^{m} \frac{1}{\lambda_k}\leq c \frac{s}{\lambda_s}
\end{equation}
with a constant $c>0.$

Suppose now that the spectral representation of $S_{\ast}$ is  $S_{\ast}=\sum_{k=1}^r \mu_k (\psi_k \otimes \psi_k),$ where $r={\rm rank}(S_{\ast})\geq 1,$ $\mu_k$ are non-zero eigenvalues of $S_{\ast}$ (possibly repeated) and $\psi_k$ are the corresponding orthonormal eigenfuctions. Denote $L:={\rm supp}(S_{\ast}).$ 
Let $\varphi $ be an arbitrary nondecreasing function such that 
$k\mapsto \frac{\varphi(k)}{k}$ is nonincreasing and 
$$
\sum_{j=1}^k \|P_L \phi_j\|^2\leq \varphi (k), k=0,1,\dots, m.
$$
We will denote by $\Psi=\Psi_{S_{\ast},W}$ the class of 
all the functions satisfying these properties. Often, it will 
be convenient to extend a function $\varphi\in \Psi$ to 
nonnegative real numbers by making it linear in each of the intervals 
$[k,k+1], k=0,1,\dots, m-1$ and setting $\varphi (u)=\varphi(m)$ for 
all $u>m.$ Such an extension will be also denoted by $\varphi.$ It is 
easy to see that the extension is a nondecreasing function in ${\mathbb R}_{+}$ and the function $u\mapsto 
\frac{\varphi (u)}{u}$ is nonincreasing. 

The following {\it coherence function} will be crucial in our analysis:
$$
\bar \varphi (k):= \bar \varphi (S_{\ast},k):=\max_{t\leq k}t \max_{j\geq t}\frac{1}{j}\sum_{i=1}^j \|P_L \phi_i\|^2,
k=1,\dots, m, \ \ \ \bar \varphi (0)=0.
$$
It is straightforward to check that $\bar \varphi \in \Psi$  
and, for all $\varphi \in \Psi,$ $\bar \varphi (k)\leq \varphi(k), k=0,\dots, m.$ 
Thus, $\bar \varphi$ is the smallest function $\varphi\in \Psi .$ 
Also, $\bar \varphi (m)=r$ since 
$\sum_{j=1}^m \|P_L \phi_j\|^2=\|P_L\|_2^2=r.$
Moreover, since $\frac{\bar\varphi(k)}{k}$ is nonincreasing,
we have 
$$\bar\varphi(k)\geq \frac{rk}{m}, k=0,\dots, m.$$

Given $t>0$, let $t_{n,m}:=t+\log(2m(\log_2(4 n^{\zeta}m^{(3/2)\zeta})+2)).$ We will assume in what follows that $m t_{n,m}\leq n$ and set
$$
\eps := 4\sqrt{\frac{t+\log (2m)}{nm}}.
$$ 

\begin{theorem}
\label{main}
There exists constants $C, C_1$ depending only on $c$ such that, for all 
$s\in \{k_0+1, \dots, m+1\}$ and all 
$\bar \eps \in [\lambda_{s}^{-1}, \lambda_{s-1}^{-1}],$\footnote{Here and in what follows, we use a convention that $\lambda_{m+1}=+\infty$ and $\lambda_{m+1}^{-1}=0.$}
with probability at least $1-e^{-t},$
\begin{equation}
\label{bound_th_main}
\|\hat S-S_{\ast}\|_{L_2(\Pi^2)}^2\leq 
C\frac{\bar \varphi(S_{\ast};s) m t_{n,m}}{n}
+\bar \eps \|W^{1/2}S_{\ast}\|_{L_2(\Pi^2)}^2
+
C_1 \max_{v\in V}\|P_L e_v\|^2 \Bigl(\frac{m t_{n,m}}{n}\Bigr)^2.
\end{equation}
\end{theorem}

{\bf Remarks}. Note that $\max_{v\in V}\|P_L e_v\|^2\leq 1.$
Thus, the last term in the righthand side of bound (\ref{bound_th_main})
is smaller than the first term, provided that 
$$\bar \varphi (S_{\ast};s)\geq \frac{mt_{n,m}}{n}.$$
Moreover, this term is much smaller under a low coherence condition 
$\max_{v\in V}\|P_L e_v\|^2\leq \frac{\nu r}{m}$ for some $\nu \geq 1$
(see conditions (\ref{coherA})). In this case,  
$$
\max_{v\in V}\|P_L e_v\|^2 \Bigl(\frac{m t_{n,m}}{n}\Bigr)^2\leq \frac{\nu r m t_{n,m}^2}{n^2}\leq \frac{\nu r t_{n,m}}{n}.
$$
Note also that Theorem \ref{main} holds in the case when $\bar \eps=0.$
In this case, $s=m$ and $\bar \varphi (S_{\ast}, m)=r,$  
so the bound of Theorem \ref{main} becomes 
\begin{equation}
\label{worse_bound}
\|\hat S-S_{\ast}\|_{L_2(\Pi^2)}^2\leq 
C\frac{r m t_{n,m}}{n},
\end{equation}
which also follows from the result of Koltchinskii, Lounici
and Tsybakov (2011) (see Theorem \ref{K-L-T} in Section 1).

The function $\bar \varphi$ involved in the statement of the theorem 
has some connection to the low coherence assumptions frequently used in the literature 
on low rank matrix completion. To be specific, suppose that, for some $\nu\geq 1,$ 
\begin{equation}
\label{low_coherence}
\sum_{j=1}^k \|P_L \phi_j\|^2 \leq \frac{\nu r k}{m},
k=1,\dots, m.
\end{equation}
Then 
$$
\bar \varphi (k)\leq \frac{\nu r k}{m}, k=1,\dots, m.
$$
A part of standard low coherence assumptions on matrix $S_{\ast}$ with respect to 
the orthonormal basis $\{\phi_k\}$ is (see (\ref{coherA}))
$$
\|P_L \phi_k\|^2 \leq \frac{\nu r}{m}, k=1,\dots, m
$$
and it implies condition (\ref{low_coherence}) that can be viewed 
as a weak version of low coherence. Under condition (\ref{low_coherence}),
the following corollary of Theorem \ref{main} holds.

\begin{corollary}
\label{cor_2}
Suppose that condition (\ref{low_coherence}) holds. 
Then, there exists a constant $C>0$ depending only on $\zeta$ such that, for all 
$s\in \{k_0+1, \dots, m+1\}$ and all 
$\bar \eps \in (\lambda_{s}^{-1}, \lambda_{s-1}^{-1}],$
with probability at least $1-e^{-t},$
$$
\|\hat S-S_{\ast}\|_{L_2(\Pi^2)}^2\leq 
C\frac{\nu r s t_{n,m}}{n}
+\bar \eps \|W^{1/2}S_{\ast}\|_{L_2(\Pi^2)}^2+
C_1 \max_{v\in V}\|P_L e_v\|^2 \Bigl(\frac{m t_{n,m}}{n}\Bigr)^2.
$$
\end{corollary}
 
Note that, if $\lambda_k \asymp k^{2\beta}$ for some $\beta>1/2,$ then 
the choice of $s$ that minimizes the bound of Corollary \ref{cor_2}
is 
$
s\asymp \left(\frac{n}{\nu r t_{n,m}}\right)^{1/(2\beta+1)}\|W^{1/2}S_{\ast}\|_{L_2(\Pi)}^{2/(2\beta+1)},
$
which, under a low coherence assumption $\max_{v\in V}\|P_L e_v\|^2\leq \frac{\nu r}{m},$ yields the bound 
\begin{equation}
\label{better_bound}
\|\hat S-S_{\ast}\|_{L_2(\Pi^2)}^2\leq 
C\biggl(\frac{\nu r t_{n,m}}{n}\biggr)^{2\beta/(2\beta+1)} \|W^{1/2}S_{\ast}\|_{L_2(\Pi)}^{2/(2\beta+1)}.
\end{equation}
The advantage of (\ref{better_bound}) comparing with (\ref{worse_bound})
(that holds for $\bar\eps=0$ and does not rely on any smoothness 
assumption on the kernel $S_{\ast}$) is due to the fact that there is no factor 
$m$ in the numerator in the right hand side of (\ref{better_bound}). 
Due to this fact, when $m$ is large enough and $\nu$ is not too large, bound (\ref{better_bound}) becomes sharper than (\ref{worse_bound}).

\section{Proofs}

{\bf Proof of Theorem \ref{main}}. Bound (\ref{bound_th_main}) will be proved 
for an arbitrary function $\varphi\in \Psi_{S_{\ast},W}$ with 
$\varphi(k)=r, k\geq m$ instead of $\bar\varphi .$ It then can be applied to the function $\bar \varphi$ (which is the smallest function in $\Psi_{S_{\ast},W}$). 
We will also assume throughout the proof that $s\in \{k_0,\dots, m\}$ and 
$\bar \eps \in [\lambda_{s+1}^{-1}, \lambda_s^{-1}]$ (at the end of the proof, 
we replace $s+1\mapsto s$). 

Denote 
$
\mathcal{P}_L(A):=A-P_{L^\perp}A P_{L^\perp},\  \mathcal{P}_L^\perp(A)=P_{L^\perp}A P_{L^\perp}, A\in {\mathcal S}_V.
$
Clearly, this defines orthogonal projectors ${\mathcal P}_L, {\mathcal P}_L^{\perp}$
in the space ${\cal S}_V$ with Hilbert--Schmidt inner product.
We will use the following well known representation of subdifferential of convex function $S \mapsto \|S\|_1:$
$
\partial \|S\|_1=\left\{ {\rm sign}(S) + \mathcal{P}_L^\perp(M) : M \in \mathcal{S}_V, \| M \| \leq 1\right\},
$
where $L={\rm supp}(S)$ (see Koltchinskii (2011b), Appendix A.4 and references 
therein).
An arbitrary matrix $A \in \partial L_n(\hat{S})$ can be represented as follows:
\begin{equation}
\label{Aravno}
A=\frac{2}{m^2} \hat{S} - \frac{2}{n}\sum_{i=1}^{n} Y_i E_{X_i,X_i'}+\eps \hat{V}+2\eps_1 W \hat{S},
\end{equation}
where $\hat{V} \in \partial \| \hat{S} \|_1$ and $E_{u,v}=E_{v,u}=\frac{1}{2}(e_u \otimes e_v + e_v \otimes e_u)$. Since $\hat{S}$ is a minimizer of $L_n(S),$ there exists a matrix $A \in \partial L_n(\hat{S})$ such that $-A$ belongs to the normal cone of $\mathbb{D}$ at the point $\hat{S}$ (see Aubin and Ekeland (1984), Chap. 2, Corollary 6).  
This implies that $\langle A, \hat{S} - S_\ast \rangle \leq 0$ and, in view of
(\ref{Aravno}),
$$
2\langle \hat{S}, \hat{S}-S_\ast\rangle_{L_2(\Pi^2)}-\left\langle \frac{2}{n} \sum_{i=1}^{n} Y_i E_{X_i,X_i'},\hat{S}-S_\ast \right\rangle+\eps \langle \hat{V},\hat{S}-S_\ast \rangle + 2\eps_1 \langle W\hat{S},\hat{S}-S_\ast \rangle \leq 0
$$
It follows by a simple algebra that 
\begin{equation}\label{eq:step1}
\begin{split}
2\|\hat{S}-S_\ast \|^2_{L_2(\Pi^2)}+2\eps_1 \|W^{1/2}(\hat{S}-S_\ast) \|^2_2+\eps \langle \hat{V},\hat{S}-S_\ast \rangle \\
\leq -2\eps_1 \langle S_\ast,W(\hat{S}-S_\ast)\rangle+2\langle \Xi,\hat{S}-S_\ast \rangle,
\end{split}
\end{equation}
where 
$$\Xi:=\frac{1}{n} \sum_{j=1}^n Y_j E_{X_j,X_j'}-{\mathbb E}YE_{X,X'}.$$
Note that $\langle \Xi, S \rangle=  \frac{1}{n} \sum_{j=1}^n \left(Y_j S(X_j,X_j')-{\mathbb E}Y S(X,X')\right)$.

On the other hand, let $V_\ast \in \partial \| S_\ast \|_1$. Therefore, the representation $V_\ast={\rm sign}(S_\ast) +\mathcal{P}_L^{\perp}(M)$ holds, where $M$ is a matrix with $\|M\| \leq 1$. It follows from the trace duality property that there exists an $M$ with $\| M \| \leq 1$ such that
$$
\langle \mathcal{P}_L^{\perp}(M), \hat{S}-S_\ast \rangle=\langle M ,\mathcal{P}_L^{\perp}(\hat{S}-S_\ast) \rangle=\langle M ,\mathcal{P}_L^{\perp}(\hat{S}) \rangle=\| \mathcal{P}_L^{\perp}(\hat{S}) \|_1
$$ 
where in the first equality we used that $\mathcal{P}_L^\perp$ is a self-adjoint operator and in the second equality we used that $S_\ast$ has support $L$. Using this equation and monotonicity of subdifferentials of convex functions, we get
$$
\langle {\rm sign} (S_\ast),\hat{S}-S_\ast \rangle+ \| \mathcal{P}_L^{\perp}(\hat{S})\|_1=\langle V_\ast, \hat{S}-S \rangle \leq \langle \hat{V}, \hat{S}-S_\ast \rangle
$$
Substituting this in \eqref{eq:step1}, it is easy to get 
\begin{eqnarray}
\label{basic}
&&
2\|\hat S-S_{\ast}\|_{L_2(\Pi^2)}^2+ \eps \|\mathcal{P}_L^{\perp}(\hat S)\|_1
+2\eps_1 \|W^{1/2}(\hat S-S_{\ast})\|_2^2 
\leq 
\\
&&
\nonumber
-\eps \langle {\rm sign}(S_{\ast}),\hat S-S_{\ast}\rangle
-2\eps_1 \langle W^{1/2}S_{\ast}, W^{1/2}(\hat S-S_{\ast})\rangle
+ 2 \langle \Xi,\hat S-S_{\ast}\rangle 
\end{eqnarray}

We will bound separately each term in the right hand side. 
First note that 
\begin{eqnarray}
\label{tri_simple}
&&
\nonumber
\eps |\langle {\rm sign}(S_{\ast}),\hat S-S_{\ast}\rangle|
\leq \eps \|{\rm sign}(S_{\ast})\|_2 \|\hat S-S_{\ast}\|_2
\\
&&
=\eps \sqrt{r}m \|\hat S-S_{\ast}\|_{L_2(\Pi^2)}
\leq \frac{1}{2}r m^2\eps^2 + \frac{1}{2}\|\hat S-S_{\ast}\|_{L_2(\Pi^2)}^2.
\end{eqnarray}
We will also need a more subtle bound on $\langle {\rm sign}(S_{\ast}),\hat S-S_{\ast}\rangle,$ expressed in terms of function $\varphi.$
Note that, for all $k_0\leq s\leq m,$
\begin{eqnarray}
&&
\nonumber
\langle {\rm sign}(S_{\ast}),\hat S-S_{\ast}\rangle=
\sum_{k=1}^m \langle {\rm sign}(S_{\ast})\phi_k, (\hat S-S_{\ast})\phi_k\rangle =
\\
&&
\nonumber
\sum_{k=1}^s \langle {\rm sign}(S_{\ast})\phi_k, (\hat S-S_{\ast})\phi_k\rangle
+
\sum_{k=s+1}^m \biggl\langle \frac{{\rm sign}(S_{\ast})\phi_k}{\sqrt{\lambda_k}}, \sqrt{\lambda_k}(\hat S-S_{\ast})\phi_k\biggr\rangle,
\end{eqnarray}
which easily implies 
\begin{eqnarray}
&&
\label{odin}
|\langle {\rm sign}(S_{\ast}),\hat S-S_{\ast}\rangle|\leq 
\biggl(\sum_{k=1}^s \|{\rm sign}(S_{\ast})\phi_k\|^2\biggr)^{1/2} 
\biggl(\sum_{k=1}^s \|(\hat S-S_{\ast})\phi_k\|^2\biggr)^{1/2}+
\\
&&
\nonumber
\biggl(\sum_{k=s+1}^m 
\frac{\|{\rm sign}(S_{\ast})\phi_k\|^2}{\lambda_k}\biggr)^{1/2}
\biggl(\sum_{k=s+1}^{m}\lambda_k\|(\hat S-S_{\ast})\phi_k\|^2\biggr)^{1/2}
\leq 
\\
&&
\nonumber
\biggl(\sum_{k=1}^s \|P_L\phi_k\|^2\biggr)^{1/2}\|\hat S-S_{\ast}\|_2
+
\biggl(\sum_{k=s+1}^m 
\frac{\|P_L\phi_k\|^2}{\lambda_k}\biggr)^{1/2}
\|W^{1/2}(\hat S-S_{\ast})\|_2.
\end{eqnarray}

We will now use the following elementary lemma.

\begin{lemma}
\label{bound_sum_s}
Let $c$ be the constant from condition (\ref{spec}). For all $s\geq k_0-1,$
$$
\sum_{k=s+1}^m 
\frac{\|P_L\phi_k\|^2}{\lambda_k}\leq 
(c+2)\frac{\varphi (s+1)}{\lambda_{s+1}}.
$$ 
\end{lemma}

{\bf Proof}. Denote $F_s:=\sum_{k=1}^s \|P_L\phi_k\|^2, s=1,\dots, m.$
Then, using the properties of function $\varphi\in \Psi,$ we get 
\begin{eqnarray}
\label{dva}
&&
\nonumber
\sum_{k=s+1}^m \frac{\|P_L\phi_k\|^2}{\lambda_k}=
\sum_{k=s+1}^{m-1} F_k \biggl(\frac{1}{\lambda_k}-\frac{1}{\lambda_{k+1}}\biggr)
+\frac{F_m}{\lambda_m}-\frac{F_s}{\lambda_{s+1}}\leq 
\\
&&
\nonumber
\sum_{k=s+1}^{m-1} 
\varphi(k)\biggl(\frac{1}{\lambda_k}-\frac{1}{\lambda_{k+1}}\biggr)
+\frac{\varphi(m)}{\lambda_m} \leq 
\frac{\varphi(s+1)}{s+1}\biggl[\sum_{k=s+1}^{m-1} 
k\biggl(\frac{1}{\lambda_k}-\frac{1}{\lambda_{k+1}}\biggr)
+\frac{m}{\lambda_m}\biggr] \leq
\\
&&
\nonumber
\frac{\varphi(s+1)}{s+1}\biggl[\sum_{k=s+2}^{m} \frac{k-(k-1)}{\lambda_k}+ \frac{(s+1)}{\lambda_{s+1}}+\frac{m}{\lambda_m}\biggr]
=
\frac{\varphi(s+1)}{s+1}\biggl[\sum_{k=s+2}^{m} \frac{1}{\lambda_k}+ \frac{(s+1)}{\lambda_{s+1}}+\frac{m}{\lambda_m}\biggr].
\end{eqnarray}
Using the assumptions on the spectrum of $W$ (in particular, condition (\ref{spec})),
we conclude that 
$$
\sum_{k=s+1}^m \frac{\|P_L\phi_k\|^2}{\lambda_k}\leq 
\frac{\varphi(s+1)}{s+1}\biggl[c\frac{s+1}{\lambda_{s+1}}+ \frac{(s+1)}{\lambda_{s+1}}+\frac{m}{\lambda_m}\biggr]\leq 
(c+2)\frac{\varphi(s+1)}{\lambda_{s+1}},
$$
ending the proof.

\qed

It follows from from (\ref{odin}) and the bound of Lemma \ref{bound_sum_s} that
\begin{eqnarray}
\label{tri}
&&
\nonumber
|\langle {\rm sign}(S_{\ast}),\hat S-S_{\ast}\rangle|\leq 
\sqrt{\varphi(s)}\|\hat S-S_{\ast}\|_2
+
\sqrt{(c+2)\frac{\varphi(s+1)}{\lambda_{s+1}}}
\|W^{1/2}(\hat S-S_{\ast})\|_2= 
\\
&&
m\sqrt{\varphi(s)}\|\hat S-S_{\ast}\|_{L_2(\Pi^2)}
+
m\sqrt{(c+2)\frac{\varphi(s+1)}{\lambda_{s+1}}}
\|W^{1/2}(\hat S-S_{\ast})\|_{L_2(\Pi^2)}. 
\end{eqnarray}
This implies the following bound:
\begin{eqnarray}
\label{tri'}
&&
\eps|\langle {\rm sign}(S_{\ast}),\hat S-S_{\ast}\rangle|\leq 
\\
&&
\nonumber
\varphi (s)m^2\eps^2+
\frac{1}{4}\|\hat S-S_{\ast}\|_{L_2(\Pi^2)}^2
+
(c+2)\frac{\varphi(s+1)}{\lambda_{s+1}}\frac{m^2\eps^2}{\bar \eps}
+\frac{\bar \eps}{4}\|W^{1/2}(\hat S-S_{\ast})\|_{L_2(\Pi^2)}^2,
\end{eqnarray}
where we used twice an elementary inequality $ab\leq a^2+\frac{1}{4}b^2, a,b>0.$
Since, under the assumptions of the theorem, $\bar\eps \lambda_{s+1}\geq 1,$ 
(\ref{tri'}) yields the following bound: 
\begin{eqnarray}
\label{tri''}
&&
\eps|\langle {\rm sign}(S_{\ast}),\hat S-S_{\ast}\rangle|\leq 
\\
&&
\nonumber
(c+3)\varphi (s+1)m^2\eps^2+
\frac{1}{4}\|\hat S-S_{\ast}\|_{L_2(\Pi^2)}^2
+\frac{\bar \eps}{4}\|W^{1/2}(\hat S-S_{\ast})\|_{L_2(\Pi^2)}^2.
\end{eqnarray}

To bound the second term in the right hand side of (\ref{basic}), 
note that 
\begin{equation}
\label{chetyre}
|\langle W^{1/2}S_{\ast}, W^{1/2}(\hat S-S_{\ast})\rangle|\leq 
\|W^{1/2}S_{\ast}\|_2 \|W^{1/2}(\hat S-S_{\ast})\|_2,
\end{equation}
which implies 
\begin{eqnarray}
\label{chetyre'}
&&
\nonumber
\eps_1 |\langle W^{1/2}S_{\ast}, W^{1/2}(\hat S-S_{\ast})\rangle|\leq 
\eps_1\|W^{1/2}S_{\ast}\|_2^2 + \frac{\eps_1}{4}\|W^{1/2}(\hat S-S_{\ast})\|_2^2
=
\\
&&
\bar \eps \|W^{1/2}S_{\ast}\|_{L_2(\Pi^2)}^2 + \frac{\bar \eps}{4}\|W^{1/2}(\hat S-S_{\ast})\|_{L_2(\Pi^2)}^2.
\end{eqnarray}

Finally, we bound $\langle \Xi, \hat{S}-S_\ast \rangle:$
\begin{equation}\label{pjat}
\begin{split}
| \langle \Xi, \hat{S}-S_\ast \rangle |& \leq 
|\langle \Xi, {\cal P}_L (\hat S-S_{\ast})\rangle| +
|\langle \Xi, {\cal P}_L^{\perp} (\hat S)\rangle|\\
& \leq |\langle {\cal P}_L\Xi, \hat S-S_{\ast}\rangle|+
\|\Xi\|\|{\cal P}_L^{\perp}(\hat S)\|_1.
\end{split}
\end{equation}
To bound $\|\Xi\|,$ we use a version of noncommutative Bernstein 
inequality of Ahlswede and Winter (2002) (see also Tropp (2010), Koltchinskii (2011a, 2011b, 2011c) for other versions of such inequalities). 

\begin{lemma}
\label{thm:OpBer}
Let $Z$ be a bounded random symmetric matrix with $\EE Z=0$, $\sigma_Z^2:=\|\EE Z^2\|$ and $\|Z\| \leq U$ for some $U > 0$. Let $Z_1,\dots,Z_n$ be $n$ i.i.d. copies of $Z$. Then for all $t>0$, with probability at least $1-e^t$
$$
\left\| \frac{1}{n}\sum_{i=1}^{n}{Z_i} \right\| \leq 2 \left( \sigma_Z  \sqrt{\frac{t+\log(2m)}{n}} \bigvee U \frac{t+\log(2m)}{n} \right)
$$
\end{lemma}

It is applied to i.i.d. random matrices $Z_i:=Y_i E_{X_i,X_i'}-\EE(Y_i E_{X_i X_i'}), i=1,\dots, n.$ Since $\|Z_i\|\leq 2$ and, by a simple computation, $\sigma_{Z_i}^2:=\|{\mathbb E}Z_i^2\|\leq 1/m$ (see, e.g., 
Koltchinskii (2011b), Section 9.4), 
Lemma \ref{thm:OpBer} implies that with probability at least $1-e^{-t}$
$$
\| \Xi \| = \left\|\frac{1}{n}\sum_{i=1}^n Z_i \right \| \leq 2\biggl[\sqrt{\frac{t+\log (2m)}{n m}}\bigvee
\frac{2(t+\log (2m))}{n}\biggr].
$$
Under the assumption that 
$$
\eps \geq  4\biggl[\sqrt{\frac{t+\log (2m)}{n m}}\bigvee
\frac{2(t+\log (2m))}{n}\biggr],
$$
this yields $\|\Xi\|\leq \eps/2$ and 
\begin{equation}
\label{pjat'}
|\langle \Xi, \hat{S}-S_\ast \rangle|\leq 
|\langle {\cal P}_L\Xi, \hat S-S_{\ast}\rangle|+
\frac{\eps}{2}\|{\cal P}_L^{\perp}(\hat S)\|_1.
\end{equation}
For simplicity, it is assumed that $n\geq 2 m(t+\log (2m)).$ 
In this case, one can take $\eps=4\sqrt{\frac{t+\log (2m)}{nm}},$
as it has been done in the statement of the theorem.

We have to bound $|\langle {\cal P}_L\Xi, \hat S-S_{\ast}\rangle|$
and we start with the following simple bound:
\begin{equation}
\label{eq:NaivePB}
\begin{split}
|\langle {\cal P}_L\Xi, \hat S-S_{\ast}\rangle| & \leq m \|\mathcal{P}_L \Xi\|_2 \|\hat{S}-S_\ast\|_{L_2(\Pi^2)}\\
& \leq m \sqrt{2r} \| \Xi \| \|\hat{S}-S_\ast\|_{L^2(\Pi^2)} \\
& \leq \frac{1}{2} m \eps \sqrt{2r} \|\hat{S}-S_\ast\|_{L^2(\Pi^2)} \\
& \leq \frac{1}{2}m^2 \eps^2 r+\frac{1}{4}\|\hat{S}-S_\ast\|_{L_2(\Pi^2)}^2,
\end{split}
\end{equation}
where we use the fact that ${\rm rank}({\cal P}_L \Xi)\leq 2r.$ 
Substituting \eqref{tri_simple}, \eqref{chetyre'}, \eqref{pjat'} and \eqref{eq:NaivePB} in \eqref{basic}, we easily get that
\begin{equation}
\label{simple_fin}
\|\hat{S}-S_\ast\|_{L_2(\Pi^2)}^2 \leq  \frac{3}{2} r \eps^2 m^2+ 
2\bar{\eps} \|W^{1/2}S_\ast\|_{L_2(\Pi^2)}^2.
\end{equation}
For $\bar \eps=0,$ this bound follows from the results 
of Koltchinskii, Lounici and Tsybakov (2011). However,
we need a more subtle bound expressed in terms of function 
$\varphi,$ which is akin to bound (\ref{tri''}).  
To this end, we will use the following lemma.

\begin{lemma} 
\label{lemma_1}
For $\delta>0,$ let $k(\delta)$ be the largest value of $k\leq m$ such that $\lambda_k^{-1}\geq \delta^2$ (if $\lambda_1^{-1}<\delta^2,$ we set $k(\delta)=0$).
For all $t>0,$ with probability at least $1-e^{-t},$
\begin{equation}
\nonumber
\sup_{\|M\|_2\leq \delta, \|W^{1/2}M\|_2\leq 1}|\langle {\cal P}_L\Xi ,M\rangle|\leq 
2\sqrt{(4c+8)}\sqrt{\frac{t}{nm}}\delta\sqrt{\varphi(k(\delta)+1)}
+
2\sqrt{2}\delta\max_{v\in V}\|P_L e_v\|\frac{t}{n},
\end{equation}
provided that $k(\delta)<m,$ and 
\begin{equation}
\nonumber
|\langle {\cal P}_L\Xi,M \rangle|\leq 
4\sqrt{2} \delta \sqrt{\frac{r t}{n m}}+ 
2\sqrt{2}\delta\max_{v\in V}\|P_L e_v\|\frac{t}{n},
\end{equation}
provided that $k(\delta)\geq m.$
\end{lemma}

{\bf Proof}. The proof is somewhat akin to the derivation of the bounds 
on Rademacher processes in terms of Mendelson's complexities used in 
learning theory (see, e.g., Proposition 3.3 in Koltchinskii (2011b)).

Note that, 
for all symmetric $m\times m$ matrices $M$,
$$
\langle {\cal P}_L \Xi, M\rangle =
\sum_{k,j=1}^m \langle {\cal P}_L \Xi, \phi_k \otimes \phi_j\rangle 
\langle M, \phi_k \otimes \phi_j\rangle.
$$
Suppose that 
$$
\|M\|_2^2=\sum_{k,j=1}^m |\langle M, \phi_k\otimes \phi_j\rangle|^2\leq  \delta^2 
$$
and 
$$
\|W^{1/2}M\|_2^2 = 
\sum_{k,j=1}^m \lambda_k |\langle M, \phi_k\otimes \phi_j\rangle|^2\leq 1.
$$
Then, it easily follows that 
$$
\sum_{k,j=1}^m \frac{|\langle M, \phi_k\otimes \phi_j\rangle|^2}{\lambda_k^{-1}\wedge \delta^2}\leq 2,
$$
which implies
\begin{eqnarray}
\label{abab}
&&
|\langle {\cal P}_L\Xi, M\rangle|
\leq 
\\
&&
\nonumber
\biggl(\sum_{k,j=1}^m (\lambda_k^{-1}\wedge \delta^2)|\langle {\cal P}_L \Xi, \phi_k \otimes \phi_j\rangle|^2\biggr)^{1/2}\biggl(\sum_{k,j=1}^m 
\frac{|\langle M, \phi_k\otimes \phi_j\rangle|^2}{\lambda_k^{-1}\wedge \delta^2}\biggr)^{1/2}\leq 
\\
&&
\nonumber
\sqrt{2}
\biggl(\sum_{k,j=1}^m (\lambda_k^{-1}\wedge \delta^2)|\langle {\cal P}_L \Xi, \phi_k \otimes \phi_j\rangle|^2\biggr)^{1/2}. 
\end{eqnarray}
Define now the following inner product:
$$
\langle M_1,M_2 \rangle_{w} := \sum_{k,j=1}^m (\lambda_k^{-1}\wedge \delta^2)
\langle M_1, \phi_k \otimes \phi_j\rangle \langle M_2, \phi_k \otimes \phi_j\rangle
$$
and let $\|\cdot\|_w$ be the corresponding norm. We will provide an upper 
bound on 
$$
\|{\cal P}_L \Xi\|_w=
\biggl(\sum_{k,j=1}^m (\lambda_k^{-1}\wedge \delta^2)|\langle {\cal P}_L \Xi, \phi_k \otimes \phi_j\rangle|^2\biggr)^{1/2}.
$$

To this end, we use a standard Bernstein type inequality for random variables 
in a Hilbert space. It is given in the following lemma. 

\begin{lemma}
\label{thm:HiBer}
Let $\xi$ be a bounded random variable with values in a Hilbert space $\mathcal{H}$. Suppose that $\EE \xi=0,$ $\EE \|\xi \|_\mathcal{H}^2=\sigma^2$ and $\|\xi\|_\mathcal{H} \leq U$. Let $\xi_1,\dots,\xi_n$ be $n$ i.i.d. copies of $\xi_i$. Then for all $t>0$, with probability at least $1-e^t$
$$
\left\| \frac{1}{n}\sum_{i=1}^n \xi_i \right\|_\mathcal{H} \leq 2\left[ \sigma \sqrt{\frac{t}{n}} \bigvee U\frac{t}{n} \right]
$$  
\end{lemma}

Applying Lemma \ref{thm:HiBer} to the random variable 
$\xi=Y \mathcal{P}_L(E_{X,X'})-\EE Y \mathcal{P}_L(E_{X,X'}),$
we get that for all $t>0,$ with probability at least $1-e^{-t},$
\begin{eqnarray}
\label{Bernstein}
&&
\|{\cal P}_L \Xi\|_w= \biggl\|
\frac{1}{n} \sum_{j=1}^n Y_j {\cal P}_L (E_{X_j,X_j'})-{\mathbb E}
Y {\cal P}_L (E_{X,X'})
\biggr\|_w\leq 
\\
&&
\nonumber 
2\biggl[{\mathbb E}^{1/2}\|Y {\cal P}_L (E_{X,X'})\|_w^2\sqrt{\frac{t}{n}}+
\Bigl\|\|Y {\cal P}_L (E_{X,X'})\|_w\Bigr\|_{L_{\infty}}\frac{t}{n}
\biggr].
\end{eqnarray}
Using the fact that $Y\in \{-1,1\},$ we get 
\begin{eqnarray}
\label{var}
&&
{\mathbb E}\|Y {\cal P}_L (E_{X,X'})\|_w^2
=
{\mathbb E}\|{\cal P}_L (E_{X,X'})\|_w^2=
\\
&&
\nonumber
{\mathbb E}\sum_{k,j=1}^m (\lambda_k^{-1}\wedge \delta^2)|\langle {\cal P}_L (E_{X,X'}), \phi_k \otimes \phi_j\rangle|^2
=
\sum_{k,j=1}^m (\lambda_k^{-1}\wedge \delta^2)
{\mathbb E}|\langle E_{X,X'}, {\cal P}_L(\phi_k \otimes \phi_j)\rangle|^2= 
\\
&&
\nonumber
\sum_{k,j=1}^m (\lambda_k^{-1}\wedge \delta^2)
m^{-2}\sum_{u,v\in V}|\langle E_{u,v}, {\cal P}_L(\phi_k \otimes \phi_j)\rangle|^2\leq 
\\
&&
\nonumber
m^{-2}\sum_{k,j=1}^m (\lambda_k^{-1}\wedge \delta^2)
\|{\cal P}_L(\phi_k \otimes \phi_j)\|_2^2
\leq 
2m^{-2}\sum_{k,j=1}^m (\lambda_k^{-1}\wedge \delta^2)
(\|P_L\phi_k\|^2+\|P_L\phi_j\|^2)=
\\
&&
\nonumber
2m^{-1}\sum_{k=1}^m (\lambda_k^{-1}\wedge \delta^2)\|P_L\phi_k\|^2
+
2m^{-2}\sum_{k=1}^m (\lambda_k^{-1}\wedge \delta^2)\sum_{j=1}^m\|P_L\phi_j\|^2=
\\
&&
\nonumber
2m^{-1}\sum_{k=1}^m (\lambda_k^{-1}\wedge \delta^2)\|P_L\phi_k\|^2
+
2m^{-2}\sum_{k=1}^m (\lambda_k^{-1}\wedge \delta^2)\|P_L\|_2^2=
\\
&&
\nonumber
2m^{-1}\sum_{k=1}^m (\lambda_k^{-1}\wedge \delta^2)\|P_L\phi_k\|^2
+
2m^{-2}r\sum_{k=1}^m (\lambda_k^{-1}\wedge \delta^2).
\end{eqnarray}
To bound ${\mathbb E}\|Y {\cal P}_L (E_{X,X'})\|_w^2$
further,
note that 
\begin{equation}
\label{hi-hi}
\sum_{k=1}^m (\lambda_k^{-1}\wedge \delta^2)\|P_L\phi_k\|^2
\leq 
\delta^2\sum_{k\leq k(\delta)} \|P_L\phi_k\|^2
+ 
\sum_{k>k(\delta)}\lambda_k^{-1} \|P_L\phi_k\|^2.
\end{equation}

Assuming that $1\leq k(\delta)\leq m-1,$ using the bound of Lemma \ref{bound_sum_s},
the fact that $\lambda_{k(\delta)+1}^{-1}< \delta^2$ and the monotonicity 
of function $\varphi,$ 
we get 
from (\ref{hi-hi}) that 
\begin{eqnarray}
\label{au}
&&
\nonumber
\sum_{k=1}^m (\lambda_k^{-1}\wedge \delta^2)\|P_L\phi_k\|^2
\leq 
\delta^2\varphi (k(\delta))+(c+2)\frac{\varphi (k(\delta)+1)}{\lambda_{k(\delta)+1}}
\leq  
\\
&&
\delta^2\varphi (k(\delta))+(c+2)\delta^2\varphi (k(\delta)+1)
\leq 
(c+3)\delta^2\varphi (k(\delta)+1).
\end{eqnarray}
It is easy to check that (\ref{au}) holds also for $k(\delta)=0$ and $k(\delta)=m$
(in the last case, $\varphi (k(\delta)+1)=r$). 
We also have 
$$
\sum_{k=1}^m (\lambda_k^{-1}\wedge \delta^2)
\leq 
\sum_{k\leq k(\delta)} \delta^2
+ 
\sum_{k>k(\delta)}\lambda_k^{-1},
$$
which, in view of condition (\ref{spec}), implies 
\begin{equation}
\label{ua}
\sum_{k=1}^m (\lambda_k^{-1}\wedge \delta^2)
\leq \delta^2 k(\delta)+ c\frac{k(\delta)+1}{\lambda_{k(\delta)}+1}
\leq (c+1)\delta^2 (k(\delta)+1).
\end{equation}
Using bounds (\ref{var}),  (\ref{au}) and (\ref{ua}), we get, 
under the condition that $k(\delta)<m,$
\begin{eqnarray}
\label{varA}
&&
{\mathbb E}\|Y {\cal P}_L (E_{X,X'})\|_w^2
\leq 
\\
&&
\nonumber
2m^{-1}(c+3)\delta^2\varphi(k(\delta)+1)+
2m^{-2}r(c+1)\delta^2(k(\delta)+1) \leq 
\\
&&
\nonumber
2m^{-1}(c+3)\delta^2\varphi(k(\delta)+1)+
2m^{-2}r(c+1)\delta^2\frac{k(\delta)+1}{\varphi(k(\delta)+1)}\varphi(k(\delta)+1) \leq 
\\
&&
\nonumber
2m^{-1}(c+3)\delta^2\varphi(k(\delta)+1)+
2m^{-2}r(c+1)\delta^2\frac{m}{\varphi(m)}\varphi(k(\delta)+1)=
\\
&&
\nonumber
(4c+8)m^{-1}\delta^2\varphi(k(\delta)+1).
\end{eqnarray}
In the case when 
$k(\delta)\geq m,$ it is easy to show that 
\begin{equation}
\label{var_L_inf}
{\mathbb E}\|Y {\cal P}_L (E_{X,X'})\|_w^2\leq 4m^{-1}\delta^2 r.
\end{equation}

We can also bound $\Bigl\|\|Y {\cal P}_L (E_{X,X'})\|_w\Bigr\|_{L_{\infty}}^2$
as follows:
\begin{eqnarray}
\label{L_inf}
&&
\Bigl\|\|Y {\cal P}_L (E_{X,X'})\|_w\Bigr\|_{L_{\infty}}^2
=
\Bigl\|\|{\cal P}_L (E_{X,X'})\|_w\Bigr\|_{L_{\infty}}^2=
\\
&&
\nonumber
\biggl\|\sum_{k,j=1}^m (\lambda_k^{-1}\wedge \delta^2)|\langle {\cal P}_L (E_{X,X'}), \phi_k \otimes \phi_j\rangle|^2
\biggr\|_{L_{\infty}}\leq 
\\
&&
\nonumber
\max_{1\leq k\leq m}(\lambda_k^{-1}\wedge \delta^2)
\max_{u,v\in V}\sum_{k,j=1}^m 
|\langle {\cal P}_L E_{u,v}, \phi_k \otimes \phi_j\rangle|^2\leq 
\\
&&
\nonumber
\max_{1\leq k\leq m}(\lambda_k^{-1}\wedge \delta^2)
\max_{u,v\in V}\|{\cal P}_L E_{u,v}\|_2^2 
\leq 
\delta^2
\max_{u,v\in V}\|{\cal P}_L(e_u\otimes e_v)\|_2^2\leq 
2\delta^2 \max_{v\in V}\|P_L e_v\|^2. 
\end{eqnarray}

If $k(\delta)<m,$ it follows from (\ref{abab}), (\ref{Bernstein}), (\ref{varA}) and (\ref{L_inf}) that with probability 
at least $1-e^{-t},$ for all symmetric matrices $M$ with $\|M\|_2\leq \delta$
and $\|W^{1/2}M\|_2\leq 1$,  
$$
|\langle {\cal P}_L\Xi,M\rangle|\leq 
2\sqrt{(4c+8)}\sqrt{\frac{t}{nm}}\delta\sqrt{\varphi(k(\delta)+1)}
+2\sqrt{2}\delta \max_{v\in V}\|P_L e_v\|\frac{t}{n}.
$$
Alternatively, if $k(\delta)\geq m,$ we use (\ref{var_L_inf}) to get  
$$
|\langle {\cal P}_L\Xi,M\rangle|\leq 
4\delta \sqrt{\frac{rt}{n m}}
+
2\sqrt{2}\delta \max_{v\in V}\|P_L e_v\|\frac{t}{n}.
$$

\qed

It follows from Lemma \ref{lemma_1} that, for all $\delta>0,$ the following bound holds with probability at least $1-e^{-t}$
\begin{eqnarray}
\label{bound_10}
&&
\sup_{\|M\|_2\leq \delta, \|W^{1/2}M\|_2\leq 1}|\langle {\cal P}_L\Xi ,M\rangle|\leq \\
&&
\nonumber
2\sqrt{(4c+8)}\sqrt{\frac{t}{nm}}
\delta\sqrt{\varphi(k(\delta)+1)}+2\sqrt{2}\delta \max_{v\in V}\|P_L e_v\|\frac{t}{n}
\end{eqnarray}
(recall that $\varphi (k)=r$ for $k\geq m,$ so, the second bound of the lemma
can be included in the first bound). 
Moreover, the bound can be easily made uniform in $\delta\in [\delta_{-}, \delta_{+}]$ for arbitrary $\delta_{-}<{\delta_{+}}.$  
To this end, take 
$\delta_j:=\delta_{+}2^{-j},j=0,1, \dots [\log_2 (\delta_+/\delta_{-})]+1$
and use (\ref{bound_10}) for each $\delta=\delta_j$ with $\bar t:=t+\log ([\log_2 (\delta_+/\delta_{-})]+2)$ instead of $t.$ An application of the union bound and monotonicity 
of the left hand side and the right hand side of (\ref{bound_10}) with respect 
to $\delta$ then implies that with probability at least $1-e^{-t}$ for all $\delta \in [\delta_{-}, \delta_{+}]$
\begin{eqnarray}
\label{bound_11}
&&
\sup_{\|M\|_2\leq \delta, \|W^{1/2} M\|_2\leq 1}|\langle {\cal P}_L\Xi ,M\rangle|\leq 
\\
&&
\nonumber
C\sqrt{\frac{\bar t}{nm}}\delta\sqrt{\varphi(k(\delta)+1)}
+4\sqrt{2}\delta \max_{v\in V}\|P_L e_v\|\frac{\bar t}{n}.
\end{eqnarray}
where $C>0$ is a constant depending only on $c.$
Indeed, by the union bound, (\ref{bound_10}) holds with 
probability at least 
$$1-([\log_2 (\delta_+/\delta_{-})]+2)e^{-\bar t}=1-e^{-t}$$ 
for all $\delta=\delta_j, j=0,\dots,[\log_2 (\delta_+/\delta_{-})]+1.$
Therefore, for all $j=0,\dots,[\log_2 (\delta_+/\delta_{-})]+1$ 
and all $\delta \in (\delta_{j+1},\delta_j]$ 
\begin{eqnarray}
\label{bound_10''}
&&
\sup_{\|M\|_2\leq \delta, \|W^{1/2}M\|_2\leq 1}|\langle {\cal P}_L\Xi ,M\rangle|\leq \\
&&
\nonumber
2\sqrt{(4c+8)}\sqrt{\frac{\bar t}{nm}}
\delta_j\sqrt{\varphi(k(\delta_{j})+1)}+
2\sqrt{2}\delta_j \max_{v\in V}\|P_L e_v\|\frac{\bar t}{n}
\end{eqnarray}
(by monotonicity of the left hand side). Note that $k(\delta_j)\leq k(\delta)\leq 
k(\delta_{j+1}).$ We can now use the fact that 
$\frac{\varphi (k)}{\lambda_k}=\frac{\varphi (k)}{k}\frac{k}{\lambda_k}$
is a nonincreasing function and the condition $\lambda_{k+1}/\lambda_k\leq c$ to 
show that
\begin{eqnarray}
&&
\nonumber
\sqrt{\frac{\bar t}{nm}}
\delta_j\sqrt{\varphi(k(\delta_{j})+1)}+
\leq 
2\sqrt{\frac{\bar t}{nm}}
\delta_{j+1}\sqrt{\varphi(k(\delta_{j+1})+1)}\leq 
\\
&&
\nonumber
2\sqrt{\frac{\bar t}{nm}}
\sqrt{\frac{\varphi(k(\delta_{j+1})+1)}{\lambda_{k(\delta_{j+1})}}}\leq 
2\sqrt{c}\sqrt{\frac{\bar t}{nm}}
\sqrt{\frac{\varphi(k(\delta_{j+1})+1)}{\lambda_{k(\delta_{j+1})+1}}}
\\
&&
\nonumber
2\sqrt{c}\sqrt{\frac{\bar t}{nm}}
\sqrt{\frac{\varphi(k(\delta)+1)}{\lambda_{k(\delta)+1}}}
\leq 
2\sqrt{c}\sqrt{\frac{\bar t}{nm}}
\delta \sqrt{\varphi(k(\delta)+1)}.
\end{eqnarray}
This and bound (\ref{bound_10''}) imply that 
\begin{eqnarray}
\label{bound_10'''}
&&
\sup_{\|M\|_2\leq \delta, \|W^{1/2}M\|_2\leq 1}|\langle {\cal P}_L\Xi ,M\rangle|\leq \\
&&
\nonumber
4\sqrt{c(4c+8)}\sqrt{\frac{\bar t}{nm}}
\delta \sqrt{\varphi(k(\delta)+1)}
+
4\sqrt{2}\delta \max_{v\in V}\|P_L e_v\|\frac{\bar t}{n},
\end{eqnarray}
which proves bound (\ref{bound_11}).

Set $\delta$ as 
$$
\delta:=
\frac{\|\hat S-S_{\ast}\|_2}{\|W^{1/2}(\hat S-S_{\ast})\|_2}
=
\frac{\|\hat S-S_{\ast}\|_{L_2(\Pi^2)}}{\|W^{1/2}(\hat S-S_{\ast})\|_{L_2(\Pi^2)}}
$$
and assume for now that $\delta \in [\delta_{-},\delta_{+}].$ 
For a particular choice of $M:=\frac{\hat{S}-S_\ast}{\|W^{1/2}(\hat{S}-S_\ast)\|_2},$ 
we get from (\ref{bound_11}) that 
\begin{equation}
\label{bound_111}
|\langle {\cal P}_L\Xi , \hat S-S_{\ast}\rangle|\leq 
C\sqrt{\frac{\bar t}{nm}}
\|\hat S-S_{\ast}\|_2\sqrt{\varphi(k(\delta)+1)}
+
4\sqrt{2}\max_{v\in V}\|P_L e_v\|\frac{\bar t}{n}\|\hat S-S_{\ast}\|_2.
\end{equation}
Suppose now that $\delta^2\geq \bar \eps.$ Since, under assumptions 
of the theorem, $\bar \eps\in (\lambda_{s+1}^{-1},\lambda_{s}^{-1}],$  
this implies that $k(\delta)\leq k(\sqrt{\bar \eps})=s$ and 
\begin{eqnarray}
\label{bound_111'}
&&
\nonumber
|\langle {\cal P}_L\Xi , \hat S-S_{\ast}\rangle|\leq 
C\sqrt{\frac{\bar t}{nm}}
\|\hat S-S_{\ast}\|_2\sqrt{\varphi(s+1)}
+
4\sqrt{2}\max_{v\in V}\|P_L e_v\|\frac{\bar t}{n}\|\hat S-S_{\ast}\|_2
=
\\
&&
\nonumber
C\sqrt{\frac{m\bar t}{n}}
\|\hat S-S_{\ast}\|_{L_2(\Pi^2)}\sqrt{\varphi(s+1)}
+
4\sqrt{2}\max_{v\in V}\|P_L e_v\|\frac{m\bar t}{n}\|\hat S-S_{\ast}\|_{L_2(\Pi)}
\leq 
\\
&&
2C^2\frac{\varphi(s+1)m\bar t}{n}+
64\max_{v\in V}\|P_L e_v\|^2
\Bigl(\frac{m\bar t}{n}\Bigr)^2+
\frac{1}{4}\|\hat S-S_{\ast}\|_{L_2(\Pi^2)}^2.
\end{eqnarray}
In the case when $\delta^2 < \bar \eps,$ we have $k(\delta)\geq k(\sqrt{\bar \eps})=s.$ 
In this case, we again use the fact that 
$\frac{\varphi (k)}{\lambda_k}$
is a nonincreasing function and the condition $\lambda_{k+1}/\lambda_k\leq c$ 
to show that 
\begin{eqnarray}
&&
\nonumber
\sqrt{\frac{\bar t}{nm}}
\|\hat S-S_{\ast}\|_2\sqrt{\varphi(k(\delta)+1)}
=
\sqrt{\frac{m \bar t}{n}}
\|W^{1/2}(\hat S-S_{\ast})\|_{L_2(\Pi^2)}
\sqrt{\delta^2\varphi(k(\delta)+1)}
\leq
\\
&&
\nonumber 
\sqrt{\frac{m \bar t}{n}}
\|W^{1/2}(\hat S-S_{\ast})\|_{L_2(\Pi^2)}
\sqrt{\frac{\varphi(k(\delta)+1)}{\lambda_{k(\delta)}}}
\leq
\sqrt{c}\sqrt{\frac{m \bar t}{n}}
\|W^{1/2}(\hat S-S_{\ast})\|_{L_2(\Pi^2)}
\sqrt{\frac{\varphi(k(\delta)+1)}{\lambda_{k(\delta)+1}}}
\leq
\\
&&
\nonumber 
\sqrt{c}\sqrt{\frac{m \bar t}{n}}
\|W^{1/2}(\hat S-S_{\ast})\|_{L_2(\Pi^2)}
\sqrt{\frac{\varphi(s+1)}{\lambda_{s+1}}}
\leq
\sqrt{c}\sqrt{\frac{m \bar t}{n}}
{\sqrt{\bar \eps}}\|W^{1/2}(\hat S-S_{\ast})\|_{L_2(\Pi^2)}
\sqrt{\varphi(s+1)}.
\end{eqnarray}
This allows us to deduce from (\ref{bound_111}) that 
\begin{eqnarray}
\label{bound_111''}
&&
|\langle {\cal P}_L\Xi , \hat S-S_{\ast}\rangle|\leq 
\\
&&
\nonumber
\sqrt{c}C\sqrt{\frac{m \bar t}{n}}
{\sqrt{\bar \eps}}\|W^{1/2}(\hat S-S_{\ast})\|_{L_2(\Pi^2)}
\sqrt{\varphi(s+1)}
+
4\sqrt{2}\max_{v\in V}\|P_L e_v\|\frac{m\bar t}{n}\|\hat S-S_{\ast}\|_{L_2(\Pi)}
\leq 
\\
&&
\nonumber
c C^2\frac{\varphi(s+1)m \bar t}{n}+
\frac{1}{4}\bar \eps\|W^{1/2}(\hat S-S_{\ast})\|_{L_2(\Pi^2)}^2
+
32\max_{v\in V}\|P_L e_v\|^2
\Bigl(\frac{m\bar t}{n}\Bigr)^2+
\frac{1}{4}\|\hat S-S_{\ast}\|_{L_2(\Pi^2)}^2.
\end{eqnarray}
It follows from bounds (\ref{bound_111'}) and (\ref{bound_111''}) that
with probability at least $1-e^{-t},$ 
\begin{eqnarray}
\label{bound_111'''}
&&
|\langle {\cal P}_L\Xi , \hat S-S_{\ast}\rangle|\leq 
(2\vee c)C^2\frac{\varphi(s+1)m\bar t}{n}+
64\max_{v\in V}\|P_L e_v\|^2
\Bigl(\frac{m\bar t}{n}\Bigr)^2+
\\
&&
\nonumber
\frac{1}{4}\|\hat S-S_{\ast}\|_{L_2(\Pi^2)}^2+
\frac{1}{4}\bar \eps\|W^{1/2}(\hat S-S_{\ast})\|_{L_2(\Pi^2)}^2,
\end{eqnarray}
provided that 
\begin{equation}
\label{in_in}
\delta=
\frac{\|\hat S-S_{\ast}\|_2}{\|W^{1/2}(\hat S-S_{\ast})\|_2}
=
\frac{\|\hat S-S_{\ast}\|_{L_2(\Pi^2)}}{\|W^{1/2}(\hat S-S_{\ast})\|_{L_2(\Pi^2)}}
\in [\delta_{-}, \delta_{+}].
\end{equation}
It remains now to substitute bounds (\ref{tri''}), (\ref{chetyre'}), (\ref{pjat'})
and (\ref{bound_111'''}) in bound (\ref{basic}) to get that with some constants 
$C>0, C_1>0$ depending only on $c$  and with probability 
at least $1-2 e^{-t}$ 
\begin{equation}
\label{eq:bound_final}
\|\hat S-S_{\ast}\|_{L_2(\Pi^2)}^2 \leq 
C \frac{\varphi (s+1)m (\bar t+t_m)}{n}
+\bar \eps \|W^{1/2}S_{\ast}\|_{L_2(\Pi^2)}^2
+ C_1\max_{v\in V}\|P_L e_v\|^2
\Bigl(\frac{m\bar t}{n}\Bigr)^2,
\end{equation}
where $t_m:=t+\log (2m).$

We still have to choose the values of $\delta_{-}, \delta_{+}$ and to handle the case when 
\begin{equation}
\label{notin}
\delta=
\frac{\|\hat S-S_{\ast}\|_2}{\|W^{1/2}(\hat S-S_{\ast})\|_2}
=
\frac{\|\hat S-S_{\ast}\|_{L_2(\Pi^2)}}{\|W^{1/2}(\hat S-S_{\ast})\|_{L_2(\Pi^2)}}
\not\in [\delta_{-}, \delta_{+}].
\end{equation}
First note that, since the largest eigenvalue of $W$ is $\lambda_m$ and it is bounded from above by $m^{\zeta},$ we have  
$$
\|W^{1/2}(\hat S-S_{\ast})\|_2\leq \sqrt{\lambda_m} \|\hat S-S_{\ast}\|_2
\leq m^{\zeta/2}\|\hat S-S_{\ast}\|_2.
$$
Thus, $\delta \geq m^{-\zeta/2}.$
Next note that 
$$
\|W^{1/2}S_{\ast}\|_{L_2(\Pi^2)}^2\leq m^{-2}m^{\zeta}\|S_{\ast}\|_2^2
\leq  m^{\zeta},
$$
where we also took into account that the absolute values of the entries of $S_{\ast}$ are bounded by $1.$ 
It now follows from (\ref{simple_fin}) that, under the assumption 
$\frac{2 mt_m}{n}\leq 1,$  
\begin{eqnarray}
&&
\nonumber
\|\hat{S}-S_\ast\|_{L_2(\Pi^2)}^2 \leq  
\frac{3}{2}r m^2 \eps^2+
2\bar{\eps} m^{\zeta}
\leq 
\\
&&
\nonumber
24 r m^2 \frac{t+\log(2m)}{n m}+2\frac{m^{\zeta}}{\lambda_s}
\leq 
12 m+2m^{2\zeta}\leq 14 m^{2\zeta},
\end{eqnarray}
which holds with probability at least $1-e^{-t}.$
Therefore, as soon as $\|W^{1/2}(\hat{S}-S_\ast)\|_{L_2(\Pi^2)} \geq n^{-\zeta},$
we have $\delta \leq 4 n^{\zeta} m^{\zeta}.$ 

We will now take $\delta_{-}:= m^{-\zeta/2}, \delta_{+}:=4n^{\zeta} m^{\zeta}.$
Then, the only case when (\ref{notin}) can possibly hold is if  $\|W^{1/2}(\hat{S}-S_\ast)\|_{L_2(\Pi^2)} \leq n^{-\zeta}.$ In this case, we can set $$\delta:=n^{\zeta} \|\hat{S}-S_\ast\|_{L_2(\Pi^2)}\in [\delta_{-},\delta_{+}]$$ 
and follow the proof of bound (\ref{bound_111'''}) replacing throughout the argument $\|W^{1/2}(\hat S-S_{\ast})\|_{L_2(\Pi^2)}$ with $n^{-\zeta}.$ This yields 
\begin{eqnarray}
\label{konec}
&&
|\langle {\cal P}_L\Xi , \hat S-S_{\ast}\rangle|\leq 
\\
&&
\nonumber
(2\vee c)C^2\frac{\varphi(s+1)m\bar t}{n}+
64\max_{v\in V}\|P_L e_v\|^2
\Bigl(\frac{m\bar t}{n}\Bigr)^2+
\frac{1}{4}\|\hat S-S_{\ast}\|_{L_2(\Pi^2)}^2+\frac{1}{4}\bar \eps n^{-2\zeta}.
\end{eqnarray}
Bound (\ref{konec}) can be now used instead of (\ref{bound_111'''}) 
to prove that 
\begin{equation}
\label{eq:bound_final''}
\|\hat S-S_{\ast}\|_{L_2(\Pi^2)}^2 \leq 
C \frac{\varphi (s+1)m (\bar t+t_m)}{n}
+\bar \eps \|W^{1/2}S_{\ast}\|_{L_2(\Pi^2)}^2 +
C_1\max_{v\in V}\|P_L e_v\|^2
\Bigl(\frac{m\bar t}{n}\Bigr)^2+
\bar \eps n^{-2\zeta}
\end{equation}
with some constants $C, C_1>0$ depending only on $c.$

Clearly, we can assume that $C_1\geq 1$ and $\bar t\geq 1.$
Since $m\leq n^2$ (recall that we even assumed that $mt_{n,m}\leq 1$),  
$\zeta\geq 1,$ $\max_{v\in V}\|P_L e_v\|^2\geq \frac{r}{m}$\footnote{Recall that $r=\|P_L\|_2^2=\sum_{v\in V}\|P_L e_v\|^2.$} 
and $\bar \eps \leq \lambda_{k_0}^{-1}\leq m^{\zeta},$ it is easy to check that 
$$
C_1\max_{v\in V}\|P_L e_v\|^2
\Bigl(\frac{m\bar t}{n}\Bigr)^2
\geq \frac{m}{n^2}\geq \frac{m^{\zeta}}{n^{2\zeta}}\geq
\bar \eps n^{-2\zeta}.
$$
Thus, the last term of bound (\ref{eq:bound_final''}) can be dropped 
(with a proper adjustment of constant $C_1$).

Note also that with our choice of $\delta_{-}, \delta_{+}$
$$
\bar t = t+\log(\log_2(\delta_{+}/\delta_{-}+2)\leq 
t+\log(\log_2(4 n^{\zeta}m^{(3/2)\zeta})+2)
$$
and $\bar t+ t_m \leq 2t_{n,m}.$ 
It is now easy to conclude that, with some constants $C, C_1$ depending 
only on $c$ and with probability at least $1-3 e^{-t}$ 
\begin{equation}
\label{eq:bound_final_A}
\|\hat S-S_{\ast}\|_{L_2(\Pi^2)}^2 \leq 
C \frac{\varphi (s+1)m t_{n,m}}{n}
+\bar \eps \|W^{1/2}S_{\ast}\|_{L_2(\Pi^2)}^2 +
C_1\max_{v\in V}\|P_L e_v\|^2
\Bigl(\frac{m\bar t}{n}\Bigr)^2.
\end{equation}
The probability bound $1-3e^{-t}$ can be rewritten as $1-e^{-t}$
by changing the value of constants $C, C_1.$ Also, by changing the notation 
$s+1\mapsto s,$ bound (\ref{eq:bound_final_A}) yields (\ref{bound_th_main}). 
This completes the proof of the theorem.

\qed




\begin{thebibliography}{99}


\bibitem{Ahlswede} Ahlswede, R. and Winter, A. (2002) Strong converse for 
identification via quantum channels. \it IEEE Transactions on Information 
Theory, \rm 48, 3, pp. 569--679.

\bibitem{Aubin} Aubin, J.-P. and Ekeland, I. (1984) Applied Nonlinear Analysis. 
J. Wiley\&Sons, New York.





\bibitem{Candes_Recht} Candes, E. and Recht, B. (2009) Exact matrix 
completion via convex optimization. \it Foundations of Computational 
Mathematics, \rm 9(6), 717--772.

\bibitem{Candes_Tao} Candes, E. and Tao, T. (2010) The power of 
convex relaxation: Near-optimal matrix completion. 
{\it IEEE Transactions on Information Theory}, 56, 2053--2080. 

\bibitem{Candes_Plan} Candes, E. and Plan, Y. (2011) Tight Oracle 
Bounds for Low-Rank Matrix Recovery from a Minimal Number of 
Random Measurements. {\it IEEE Transactions on Information Theory},
57(4), 2342--2359.









\bibitem{Gross-2} Gross, D. (2011) Recovering Low-Rank Matrices From 
Few Coefficients in Any Basis. {\it IEEE Transactions on Information Theory},
57, 3, 1548--1566. 



\bibitem{Ko_25} Koltchinskii, V. (2011a) Von Neumann Entropy Penalization 
and Low Rank Matrix Estimation. \it Annals of Statistics, \rm 39, 6, 2936--2973.


\bibitem{Ko_27} Koltchinskii, V. (2011b) Oracle Inequalities 
in Empirical Risk Minimization and Sparse Recovery Problems,
\it Ecole d'ete de Probabilit\'es de Saint-Flour 2008, \rm
Lecture Notes in Mathematics, Springer. 

\bibitem{Ko_28} Koltchinskii, V. (2011c) A remark on low rank matrix 
recovery and noncommutative Bernstein type inequalities. Preprint.



\bibitem{Ko_26} Koltchinskii, V., Lounici, K. and Tsybakov, A. (2011)
Nuclear norm penalization and optimal rates for noisy matrix completion.
\it Annals of Statistics, \rm 39, 5, 2302--2329.








\bibitem{Negahban} Negahban, S. and Wainwright, M.J. (2010)
Restricted strong convexity and weighted matrix completion 
with noise. Preprint. 


\bibitem{Recht_1} Recht, B., Fazel, M. and Parrilo, P. (2010) Guaranteed 
minimum rank solutions of matrix equations via nuclear norm minimization.
{\it SIAM Review}, 52, 3, 471--501. 



\bibitem{Rohde} Rohde, A. and Tsybakov, A. (2011) Estimation of high-dimensional low rank matrices. {\it Annals of Statistics}, 39, 2, 887--930.








\bibitem{Tropp} Tropp, J.A. (2010) User-friendly tail bounds for sums of 
random matrices. {\it Foundations of Computational Mathematics}, to appear.



\end{thebibliography}
\end{document}